\documentclass[final]{siamart171218} 

\usepackage[utf8]{inputenc} 
\usepackage{amsmath}
\usepackage{amssymb}
\usepackage{color}
\usepackage{subfig}
\usepackage{graphicx}
\usepackage{enumitem}
\usepackage{booktabs}
\usepackage{array}
\definecolor{hotpink}{rgb}{0.9,0,0.5}
\hypersetup{colorlinks,urlcolor=blue,citecolor=hotpink,linkcolor=blue}

\usepackage{listings}

\title{The Power of Bidiagonal Matrices\thanks{Version of November 6, 2023.
This paper is based on the Hans Schneider Prize talk given at the 
25th International Linear Algebra Society (ILAS),
Conference, Madrid, June 12--16, 2023.
This research was supported by the Royal Society.}}
\author{Nicholas J. Higham%
        \thanks{%
                Department of Mathematics,
                University of Manchester,
                Manchester, M13 9PL, UK
                (\texttt{nick.higham@manchester.ac.uk}).
               }
}

\usepackage[trueslanted]{newtxtext}  
\usepackage[smallerops]{newtxmath}
\DeclareSymbolFont{CMlargesymbols}{OMX}{cmex}{m}{n}
\DeclareMathDelimiter{(}{\mathopen} {operators}{"28}{CMlargesymbols}{"00}
\DeclareMathDelimiter{)}{\mathclose}{operators}{"29}{CMlargesymbols}{"01}
\DeclareMathAlphabet\mathcal{OMS}{cmsy}{m}{n}
\SetMathAlphabet\mathcal{bold}{OMS}{cmsy}{b}{n}

\def\shat{\widehat{s}}
\def\csp{checkerboard sign pattern}

\def\smatrix#1{\bigl[\begin{smallmatrix}#1\end{smallmatrix}\bigr]}

\def\eu{\ensuremath{\mathrm{e}}}

\def\rhs{right-hand side}

\def\KMS{Kac--Murdock--Szeg\"o}
\def\DS{\displaystyle}

\newcommand{\Wlogy}{Without loss of generality}
\newcommand{\wlogy}{without loss of generality}

\def\berr{backward error}
\def\ferr{forward error}
\def\crbe{componentwise relative backward error}

\def\py{polynomial}

\def\eval{eigenvalue}
\def\evals{eigenvalues}

\def\evecs{eigenvectors}

\def\sval{singular value}
\def\svals{singular values}

\def\tri{triangular}
\def\cmp{componentwise}

\def\pert{perturbation}

\def\resp{respectively}

\def\TN{totally nonnegative}
\def\TP{totally positive}

\def\sval{singular value}
\def\svals{singular values}

\def\resp{respectively}

\def\eps{\epsilon}
\def\l{\lambda}

\DeclareMathOperator{\fl}{\operatorname{f\kern.2ptl}} 
\DeclareMathOperator{\cond}{\operatorname{cond}}
\DeclareMathOperator{\diag}{\operatorname{diag}}

\DeclareMathOperator{\sign}{\operatorname{sign}}

\def\normt#1{\|#1\|_2}
\def\normi#1{\|#1\|_1}

\def\Sig{\Sigma} 
\def\s{\sigma}

\makeatletter
\def\mycases#1{\left\{\,\vcenter{\normalbaselines\m@th
    \ialign{$##\hfil$&\quad{##}\hfil\crcr#1\crcr}}\right.}
\makeatother

\def\spd{symmetric positive definite}

\def\fact{factorization}
\def\LUf{LU \fact}

\def\k{\kappa}
\def\ktwo{\kappa_2}

\def\kinf{\kappa_{\infty}}
\def\cn{condition number}

\def\gu{\gamma^{}}

\def\berr{backward error}

\def\fpa{floating-point arithmetic}

\def\rea{rounding error analysis}
\def\mult{multiplication}

\def\xhat{\widehat{x}}
\def\yhat{\widehat{y}}

\def\chat{\widehat{c}}

\def\d{\delta}
\def\dA{\Delta A}
\def\dB{\Delta B}
\def\db{\Delta b}

\def\gu{\gamma^{}}
\def\R{\mathbb{R}}
\def\C{\mathbb{C}}
\def\Z{\mathbb{Z}}
\def\mbyn{m \times n}
\def\mbym{m \times m}
\def\nbyn{n \times n}

\def\fact{factorization}

\def\normi#1{\|#1\|_1}
\def\normt#1{\|#1\|_2}

\def\normo#1{\|#1\|_{\infty}}
\def\normoq#1{\|\,#1\,\|_{\infty}}

\def\Alg{Algorithm}
\def\alg{algorithm}
\def\a{\alpha}
\def\b{\beta}
\def\th{\theta}

 \makeatletter
 \def\mymatrix#1{\null\,\vcenter{\normalbaselines\m@th
     \ialign{\hfil$##$\hfil&&\quad\hfil$##$\hfil\crcr
       \mathstrut\crcr\noalign{\kern-\baselineskip}
       #1\crcr\mathstrut\crcr\noalign{\kern-\baselineskip}}}\,}
 \makeatother
 \def\Bmatrix#1{\left[ \mymatrix{#1} \right]}

\mathchardef\Gamma="7100 \mathchardef\Delta="7101
\mathchardef\Theta="7102 \mathchardef\Lambda="7103
\mathchardef\Xi="7104 \mathchardef\Pi="7105 \mathchardef\Sigma="7106
\mathchardef\Upsilon="7107 \mathchardef\Phi="7108
\mathchardef\Psi="7109 \mathchardef\Omega="710A
\definecolor{gray}{rgb}{0.5,0.5,0.5}
\definecolor{mauve}{rgb}{0.58,0,0.82}
\definecolor{lightgrey}{rgb}{0.9,0.9,0.9}
\definecolor{darkgreen}{rgb}{0,0.6,0}

\mathcode`@="8000 
{\catcode`\@=\active\gdef@{\mkern1mu}}

\renewtheorem{algorithm}[theorem]{Algorithm}
\newcounter{mylineno}
\makeatletter
\let\oldtabcr\@tabcr
\def\nonumberbreak{\oldtabcr\hspace{3.5pt}}
\def\mynewline{\refstepcounter{mylineno}%
               \llap{\footnotesize\arabic{mylineno}\hspace{5pt}}%
              }

\gdef\@tabcr{\@stopline \@ifstar{\penalty%
           \@M \@xtabcr}\@xtabcr\mynewline}

\newenvironment{code}{%
                        \mathcode`\:="603A  
                        \def\colon{\mathchar"303A}
                        \setcounter{mylineno}{0}
                        \par
                        \upshape
                        \begin{list} 
                           {} {\leftmargin = 1cm}
                        \item[]
                        \begin{tabbing}

                           \hspace*{.3in} \= \hspace*{.3in} \=
                           \hspace*{.3in} \= \hspace*{.3in} \= \kill
                           \mynewline
                       }{\end{tabbing}\end{list}}
\makeatother

\begin{document}
\maketitle

\begin{abstract}
Bidiagonal matrices are widespread in numerical linear algebra, 
not least because of their use in the standard algorithm
for computing the singular value decomposition and their appearance 
as LU factors of tridiagonal matrices.
We show that bidiagonal matrices have a number of interesting properties 
that make them powerful tools in a variety of problems,
especially when they are multiplied together.
We show that the inverse of a product of bidiagonal matrices is insensitive
to small componentwise relative perturbations in the factors if the factors
or their inverses are nonnegative.
We derive componentwise rounding error bounds for the solution of a 
linear system $Ax = b$,
where $A$ or $A^{-1}$ is a product
$B_1 B_2\dots B_k$ of bidiagonal matrices,
showing that strong results are obtained when the 
$B_i$ are nonnegative or have a checkerboard sign pattern.
We show that given the \fact\ of an $n\times n$  totally nonnegative matrix $A$
into the product of bidiagonal matrices,
$\normo{A^{-1}}$ can be computed in $O(n^2)$ flops
and that in floating-point arithmetic the computed result has small relative
error, no matter how large $\normo{A^{-1}}$ is.
We also show how factorizations involving bidiagonal matrices
of some special matrices,
such as the Frank matrix and the Kac--Murdock--Szeg\"o matrix, 
yield simple proofs
of the total nonnegativity and other properties of these matrices.
\end{abstract}


\section{Introduction}

Bidiagonal matrices
\begin{equation}\notag
 B = 
     \begin{bmatrix} 
        b_{11} & b_{12} &        &            \\
               & b_{22} & \ddots &            \\
               &        & \ddots & b_{n-1,n}  \\
               &        &        & b_{nn} 
     \end{bmatrix} \in\C^{\nbyn}
\end{equation}
have $2n-1$ parameters, appearing on two diagonals.
Despite their simplicity,
bidiagonal matrices are powerful tools in a variety of problems, especially
when they are 
multiplied together.
Their properties and uses have been explained by various authors,
but the full range of them 
may be underappreciated.
Indeed, in the 1139-page book \emph{Matrix Mathematics} \cite{bern09} 
the word ``bidiagonal'' appears on only one page
and bidiagonal matrices appear little in the 
\emph{Handbook of Linear Algebra}
\cite{hogb14} apart from in the chapter by Fallat~\cite{fall14}.

The purpose of this work is to show the utility of bidiagonal matrices,
and in particular to show how factorizations of matrices into bidiagonal
factors can be exploited.
Our main contributions are as follows, where $A = B_1B_2\dots B_k$
with each $B_i$ either upper bidiagonal or lower bidiagonal.
\begin{itemize}

\item 
We show that small \cmp\ \pert s in the $B_i$ produce small \cmp\ \pert s in
$A^{-1}$ if the $B_i$ or the $B_i^{-1}$ are nonnegative
(Theorem~\ref{thm.bidiag-prod-inv-pert}).

\item
We show that the \cn\ $\kinf(A) = \normo{A}\normo{A^{-1}}$ can be computed in
$O(kn)$ flops when 
the $B_i$ are nonnegative or have a \csp,
without explicitly forming $A$ (section~\ref{sec.cond-comp}).

\item
We give a unified derivation of \berr\ bounds and \ferr\ bounds for the
computed solution of $Ax = b$ when $A$ or $A^{-1}$ is a product of
bidiagonal matrices and the system is solved using the factors
(section~\ref{sec.linear-systems}).

\item
We show that for a \TN\ $\nbyn$ matrix $A$, $\kinf(A)$ can be computed 
in $O(n^2)$ flops, given a \fact\ of $A$ into a product of bidiagonal
matrices
and that the computed solution is highly accurate
(\Alg~\ref{alg.fastTNcond}).

\item
We explore functions of bidiagonal matrices and show that the exponential
of a \TN\ bidiagonal matrix is \TN.

\item
We give new observations on how \fact s involving bidiagonal matrices
can help us to understand properties of some well-known matrices
(section~\ref{sec.expl-fact}).

\end{itemize}

Bidiagonal matrices arise in some classical contexts in numerical linear
algebra, which we briefly summarize as they will not be the focus of our attention.

\emph{Computing the singular value decomposition (SVD)}.
The first step of the Golub--Reinsch \alg\ for computing the SVD is a
two-sided reduction by Householder transformations to upper bidiagonal form
$B$, 
as proposed by Golub and Kahan \cite{goka65}.
The SVD of $B$ is then computed by the QR \alg\
implicitly applied to $B^*B$, and this 
can be done in a way that guarantees high relative
accuracy in all the computed \svals\ of $B$ \cite{deka90}.

\emph{\LUf\ of tridiagonal matrices.}
If $A\in\C^{\nbyn}$ is tridiagonal and has an \LUf\ $A = LU$ then $L$ is
unit lower bidiagonal and $U$ is upper bidiagonal.

\emph{Lanczos bidiagonalization.}
For large, sparse matrices the solution to a linear system or the 
least squares solution to an overdetermined system can be computed using a method 
based on unitary reduction to bidiagonal form 
by the Lanczos process \cite[sec.~7.6]{bjor96}, 
\cite{goka65}, \cite{pasa82}.

In \pert\ and rounding error analyses products of terms of the form 
$1 + \d_i$ arise.
Their distance from $1$ will be bounded using the following result
\cite[Lem.~3.1]{high:ASNA2}.

\begin{lemma}\label{lem:gammadet}
If $|\delta_i| \le \d$ and $\rho_i = \pm 1$ for $i = 1\colon n$, 
and $n\d < 1$, then
\begin{equation}\label{thetan}
  \prod_{i=1}^n(1+\delta_i)^{\rho_i} = 1 + \theta_n, 
       \quad |\theta_n| \le \frac{n\d}{1-n\d}. 
\end{equation}   
\end{lemma}

We also need a \cmp\ bound for \pert s in a matrix product
\cite[Lem.~3.8]{high:ASNA2}.
Here and throughout,
$|A| = (|a_{ij}|)$ and inequalities between matrices hold \cmp.

\begin{lemma}\label{lem.matprod-comp}
If $X_j+\Delta X_j \in \C^{\nbyn}$ satisfies
$|\Delta X_j| \le \d_j |X_j|$ for $j = 1\colon m$ then
$$
    \biggl| \prod_{j=1}^m (X_j + \Delta X_j) - \prod_{j=1}^m X_j \biggr|
    \le
    \Biggl( \prod_{j=1}^m (1+\d_j) -1 \Biggr)
    \prod_{j=1}^m |X_j|.
$$
\end{lemma}

We use the standard model of \fpa\ \cite[sec.~2.2]{high:ASNA2}
and denote by $u$ the unit roundoff.
We need the constant, for $nu < 1$,
\begin{equation}\notag
      \gu_n = \frac{nu}{1-nu}.
\end{equation}

We will make use of the one-parameter bidiagonal matrix
\begin{equation}\label{Tnth}
  T_n(\th) = \Bmatrix{  1 &  \th &      &        & \cr
                          & 1    &  \th &        & \cr
                          &      & 1    & \ddots & \cr
                          &      &      & \ddots &  \th\cr
                          &      &      &        & 1 \cr} \in\C^{\nbyn}.
\end{equation}

\section{Basic Properties of Bidiagonal Matrices}

First we consider the inverse of a nonsingular bidiagonal matrix.
It is instructive to look at the $4 \times 4$ case:
\begin{equation}\notag
\begin{bmatrix}
a           & x               & 0                    & 0 \\   
            & b               & y                    & 0 \\   
            &                 & c                    & z \\
            &                 &                      & d 
\end{bmatrix}^{-1}
=
\begin{bmatrix}
\frac{1}{a} & -\frac{x}{a\,b} & \frac{x\,y}{a\,b\,c} & -\frac{x\,y\,z}{a\,b\,c\,d} \\[3pt] 
            & \frac{1}{b}     & -\frac{y}{b\,c}      & \frac{y\,z}{b\,c\,d}        \\[3pt] 
            &                 & \frac{1}{c}          & -\frac{z}{c\,d}             \\[3pt]
            &                 &                      & \frac{1}{d} 
\end{bmatrix}.
\end{equation}
Notice that every element in the upper triangle
is a product of off-diagonal elements of $B$
and inverses of diagonal elements,
that the superdiagonals have alternating signs attached, 
and that there are no additions.
These properties hold for general $n$,
as the explicit form of the inverse in the following result shows.

\begin{lemma}\label{lem.bidiag-inv}
If $B\in\C^{\nbyn}$ is nonsingular and upper bidiagonal then 
\begin{equation}\label{Binv-ij}
   (B^{-1})_{ij} = \frac{1}{b_{jj}} 
                   \prod_{k=i}^{j-1} \biggl(\frac{-b_{k,k+1}}{b_{kk}}\biggr), \quad j\ge i.
\end{equation}
\end{lemma}

We will make use of the fact that when $B$ has nonnegative elements,
$B^{-1}$ has a checkerboard (alternating) sign pattern.

We introduce the comparison matrix $M(A)$ of $A\in\C^{\nbyn}$:
\begin{equation*}
  \bigl(M(A)\bigr)_{ij} =
   \begin{cases} |a_{ii}|, & i=j, \\
                -|a_{ij}|, & i\ne j.
    \end{cases}
\end{equation*}
It is easy to see that 
\begin{equation}\label{B-MB}
   |B^{-1}| = M(B)^{-1},
\end{equation}
an observation that we will need later.

Using the representation \eqref{Binv-ij} of the inverse we can bound the
effect of a \cmp\ \pert\ of $B$.
Let
\begin{equation}\label{tau}
  \tau = \frac{(2n-1)\d}{1 -(2n-1)\d}.
\end{equation}

\begin{theorem}\label{thm.bidiag-inv-pert}
If $B\in\C^{\nbyn}$ is a nonsingular bidiagonal matrix and $\dB$ is a
\pert\ satisfying $|\dB| \le \d |B|$ then 
\begin{equation}\notag
   |(B+\dB)^{-1} - B^{-1}\bigr|
      \le \tau |B^{-1}|,
\end{equation}
where $\tau$ is defined in \eqref{tau}.
\end{theorem}

\begin{proof}
Assume, \wlogy, that $B$ is upper bidiagonal.
Write $\db_{ij} = \d_{ij} b_{ij}$, where $|\d_{ij}| \le \d$.
From \eqref{Binv-ij} we obtain 
\begin{align*}
 (B+\dB)^{-1}_{ij} - (B^{-1})_{ij}
 & = \frac{1}{b_{jj}(1+\d_{jj})} 
      \prod_{k=i}^{j-1} \biggl(\frac{-b_{k,k+1}(1+\d_{k,k+1})}
                                    {b_{kk}(1+\d_{kk})}\biggr)
    - \frac{1}{b_{jj}} 
        \prod_{k=i}^{j-1} \biggl(\frac{-b_{k,k+1}}{b_{kk}}\biggr)\\
 & = (B^{-1})_{ij} \biggl( \frac{1}{1+\d_{jj}} 
      \prod_{k=i}^{j-1} \biggl(\frac{1+\d_{k,k+1}}
                                    {1+\d_{kk}} \biggr) -1 \biggr)\\
 & = (B^{-1})_{ij} \th_{2(j-i)+1},
\end{align*}
where 
$|\th_k| \le \gu_k = k\d/(1- k\d)$ by Lemma~\ref{lem:gammadet}.
\end{proof}

This result, which is essentially the same as 
\cite[Prob.~22.8]{high:ASNA2},
 shows that a \cmp\ relative \pert\ in $B$ produces a
\cmp\ relative \pert\ in $B^{-1}$ at most about $2n$ times larger: a strong
result that does not hold for \tri\ matrices in general.

We now extend this result to a product of bidiagonal matrices.
In all the products of bidiagonal matrices in this paper each matrix can be
upper bidiagonal or lower bidiagonal.

\begin{theorem}\label{thm.bidiag-prod-inv-pert}
Let $B = B_1B_2\dots B_k\in\C^{\nbyn}$,
where the $B_i$ are nonsingular bidiagonal matrices,  
and let  $B+\dB = (B_1+\dB_1)(B_2+\dB_2)\dots (B_k+\dB_k)$,
where $|\dB_i|\le \d |B_i|$ for all $i$.
Then
\begin{equation}\label{Binv-prod-pert}
   \bigl|(B+\dB)^{-1} - B^{-1}\bigr|
      \le \bigl( ( 1+ \tau )^k -1 \bigr)
        |B_k^{-1}| |B_{k-1}^{-1}| \dots |B_1^{-1}|,
\end{equation}
where $\tau$ is defined in \eqref{tau},
and if the $B_i$ or the $B_i^{-1}$ are all nonnegative
then 
\begin{equation}\label{B-prod-pert-inv}
   \bigl|(B+\dB)^{-1} - B^{-1}\bigr|
      \le \bigl( ( 1+ \tau )^k -1 \bigr) |B^{-1}|.
\end{equation}
\end{theorem}

\begin{proof}
We have 
\begin{align*}\notag
 (B+\dB)^{-1} 
  &= (B_k + \dB_k)^{-1} (B_{k-1} + \dB_{k-1})^{-1} \dots(B_1+\dB_1)^{-1}\\
  &= (B_k^{-1} + E_k) (B_{k-1}^{-1} + E_{k-1} ) \dots (B_1^{-1} + E_1),
\end{align*}
where by Theorem~\ref{thm.bidiag-inv-pert},
$|E_i| \le \tau |B_i^{-1}|$, $i = 1\colon k$.
Hence by Lemma~\ref{lem.matprod-comp},
\begin{equation}\notag 
 |(B+\dB^{-1}) - B^{-1}| 
               \le \bigl( ( 1+ \tau )^k -1 \bigr) 
                |B_k^{-1}|@ |B_{k-1}^{-1}| \dots |B_1^{-1}|.
\end{equation}
The bound \eqref{B-prod-pert-inv} 
is immediate if the $B_i^{-1}$ are all nonnegative.
If the $B_i$ are all nonnegative. then 
\eqref{B-prod-pert-inv} follows from considering the \csp\ of the
inverses; see Theorem~\ref{thm.Binv-prod} below.
\end{proof}

The bound \eqref{B-prod-pert-inv} shows that if the $B_i$ or the $B_i^{-1}$
are all nonnegative then 
\cmp\ relative \pert s in the
$B_i$ produce \cmp\ relative \pert\ in the inverse of the product
at most about a factor $2nk$ times larger.

Like the inverse, the \svals\ of a bidiagonal matrix are very well behaved
under \cmp\ \pert s.
Let $\s_i(B)$ denote the $i$th largest \sval\ of $B$.

\begin{theorem}\label{thm.bidiag-sval-pert}
Let $B\in\C^{\nbyn}$ and $B+\dB$ be upper bidiagonal and suppose that 
$(B+\dB)_{ii} = \a_{2i-1}b_{ii}$ and 
$(B+\dB)_{i,i+1} = \a_{2i}b_{i,i+1}$,
where the $\a_i$ are nonzero.   Then
\begin{equation}\notag
   \frac{\s_i(B)}{\mu} \le \s_i(B+\dB) \le \mu@ \s_i(B), \quad i=1\colon n,
\end{equation}
where
\begin{equation}\notag
       \mu = \prod_{i=1}^{2n-1} \max( |\a_i|, |\a_i^{-1}| ).
\end{equation}
\end{theorem}

\begin{proof}
We can write $B + \dB = D_1BD_2$, where
\begin{equation}\notag
    D_1 = \diag \Bigl( \a_1, \frac{\a_1\a_3}{\a_2},
                             \frac{\a_1\a_3\a_5}{\a_2\a_4}, \dots \Bigr),\qquad
    D_2 = \diag \Bigl(1, \frac{\a_2}{\a_1}, \frac{\a_2\a_4}{\a_1\a_3},
                             \frac{\a_2\a_4\a_6}{\a_1\a_3\a_5}, \dots \Bigr).
\end{equation}
An extension for \svals\ of a result of 
Ostroswki for \evals\ \cite[Thm.~3.1]{eiip95} gives
\begin{equation}\notag
     \frac{\s_i(B)}{\normt{D_1^{-1}} \normt{D_2^{-1}}}
    \le  \s_i(B+\dB) \le \s_i(B) \normt{D_1} \normt{D_2}.
\end{equation}
Using 
$\normt{D_1}\normt{D_2} 
= \max_i| (D_1)_{ii} | \max_i| (D_2)_{ii} |\le \mu$ 
(taking account of cancellation in the product) and
$\normt{D_1^{-1}}\normt{D_2^{-1}} \le \mu$ gives the result.
\end{proof}

Theorem~\ref{thm.bidiag-sval-pert} is from 
Demmel and Kahan \cite[Cor.~2]{deka90}
and the proof is from 
Eisenstat and Ipsen \cite[Cor.~4.2]{eiip95}.
The theorem shows that relative \pert s of magnitude at most 
$\tau = \max_i |1-\a_i| \ll 1$ 
to the elements on the diagonal and superdiagonal of an upper
bidiagonal matrix produce relative changes of at most 
$(1-\tau)^{2n-1}-1 \approx (2n-1)\tau$ in each \sval.
This is a much stronger result than for general \pert s of a general
$\nbyn$ matrix,
where it is only the absolute changes in the \svals\ that are bounded:
$|\sigma_k(A+\dA) - \sigma_k(A)| \le \s_1(\dA) = \|\dA\|_2$, 
$k = 1\colon n$ \cite[Cor.~7.3.5]{hojo13}.

Theorem~\ref{thm.bidiag-sval-pert} does not extend to a product of
bidiagonal matrices, as the following example shows.
Let 
\begin{align*}\notag
  A &= I =  \Bmatrix{1 & x \cr 0 & 1 }
      \Bmatrix{1 & -x \cr 0 & 1 } =: B_1B_2, \\
  A+\dA &= \Bmatrix{1 & 2x\d      \cr 0 & 1 }
         = \Bmatrix{1 & x(1+\d) \cr 0 & 1 }
          \Bmatrix{1 & -x(1-\d) \cr 0 & 1 } =: (B_1+\dB_1)(B_2+\dB_2),
\end{align*}
where $\d > 0$, $x > 0$, and $x\d \gg 1$.
Here, $B_1$ and $B_2$ have undergone a \cmp\ relative change $\d$.
The \svals\ of $A$ are $\s_1 = 1$ and $\s_2 = 1$, and those of $A+\dA$ are
approximately 
$\shat_1 = 2x\d$ and 
$\shat_2 = (2x\d)^{-1}$ (since $x\d \gg 1$).
Hence the relative change in $\s_1$ is
$|\s_1-\shat_1|/\s_1 \approx 2x\d \gg 1$ and that in $\s_2$ is
$|\s_2-\shat_2|/\s_2 \approx 1 - 1/(2x\d) \approx 1$.
We conclude that 
relative changes in bidiagonal matrices 
$B_1, B_2,\dots,B_k$ can induce a much larger
relative change in the \svals\ of their product.
The situation is different for a product of \emph{nonnegative}
bidiagonal matrices $B_1,B_2,\dots B_k$:
small \cmp\ relative changes in the 
$B_i$ produce only small relative changes in the \svals\
of the product $B_1,B_2,\dots B_k$,
as shown by Koev \cite[Cor.~7.3]{koev05}.

The next result reveals some further interesting properties of the \svals\
of a bidiagonal matrix.

\begin{theorem}\label{thm.SVD-bid}
Let $B\in\C^{\nbyn}$ be bidiagonal.
\begin{enumerate}[label=\upshape(\alph*)]

\item\label{SVD-bid-a}
$|B| = D B F$, where $D$ and $F$ are unitary diagonal matrices.
Hence $B$ and $|B|$ have the same \svals.

\item\label{SVD-bid-b}
If $b_{ii}$ and $b_{i,i+1}$ are nonzero for all $i$ then the \svals\ of $B$
are distinct.
\end{enumerate}

\end{theorem}

\begin{proof}
\ref{SVD-bid-a}:
Let $D = \diag(d_i)$ and $F = \diag(f_i)$ with $f_1 = 1$.
We take 
$d_1 = \sign(b_{11})^*$,
$f_2 = \sign(d_1 b_{12})^*$,
$d_2 = \sign(b_{22}f_2)^*$,
$f_3 = \sign(d_2 b_{23})^*$,
and so on, where $\sign(z)  = z/|z|$ if $z\ne 0$ or $1$ otherwise.
Then $|B| = DBF$, 
where $D$ and $F$ have diagonal elements of modulus $1$ and so are
unitary.
Therefore 
if  $B = U\Sig V^*$ is an SVD of $B$ then 
$|B| = (DU)\Sig (V^*F)$ is an SVD of $|B|$.

\ref{SVD-bid-b}:
The \svals\ of $B$ are the square roots of the \evals\ of 
$T = |B|^*|B|$, by \ref{SVD-bid-a}.
The matrix $T$ is symmetric tridiagonal with positive superdiagonal and subdiagonal
elements, so the \evals\ of $T$ are distinct
\cite[Lem.~7.7.1]{parl98},
and hence so are the \svals\ of $B$.
\end{proof}

It is interesting to note that the SVD codes in 
both LINPACK \cite{dbms79} and LAPACK \cite{lug99} reduce 
$A\in\C^{\mbyn}$ to a real bidiagonal matrix, 
so that the QR iteration can be carried out in real arithmetic,
but they do so in different ways.
LINPACK reduces $A$ to bidiagonal form by Householder transformations
and then explicitly carries out the diagonal scaling given in part 
\ref{SVD-bid-a} of Theorem~\ref{thm.SVD-bid}.
LAPACK reduces $A$ to bidiagonal form
using elementary unitary matrices of the form
$P = I - \rho vv^*$ with generally nonreal $\rho$
that are chosen so that the reduced bidiagonal matrix is real
\cite{leho96}.

\section{The Condition Number of a Matrix Product}\label{sec.cond-comp}

Suppose a matrix $X\in\C^{\nbyn}$ is given in factored form
$X = A_1A_2\dots A_k$, where $A_i\in\C^{\nbyn}$ for all $i$,
and that we wish to compute or estimate 
the \cn\ $\kinf(X) = \normo{X} \normo{X^{-1}}$ 
without explicitly forming $X$.
Initially we will make no assumptions about the $A_i$, 
but later we will specialize to bidiagonal $A_i$.
For dense matrices the cost of forming $X$ is 
$2(k-1)n^3$ flops,
whereas we would like to compute or estimate $\kinf(X)$ at the cost of a
few matrix--vector products with $X$, that is, 
in a small multiple of $2(k-1)n^2$ flops.

The \cn\ estimation problem is well studied \cite[Chap.~15]{high:ASNA2}.
Here we focus on the problem of \emph{exactly} computing the \cn.
Recall that the $\infty$-norm satisfies
\begin{equation}\notag
 \normo{X} = \normoq{|X|} = \normoq{|X|e}, 
\end{equation}
where $e = [1,1,\dots,1]^T$.

In general we cannot compute 
$\normo{ A_1A_2\dots A_k}$ without forming the matrix product.
However, if the equality
\begin{equation}\label{Aprod-eq}
 |A_1A_2\dots A_k| = |A_1|@|A_2|\dots |A_k| 
\end{equation}
holds then 
\begin{equation}\label{prod-bound}
 \normo{ A_1A_2\dots A_k } 
    = \normoq{ |A_1|@|A_2|\dots |A_k| } 
    = \normoq{ |A_1|@|A_2|\dots |A_k|e}
\end{equation}
and we can evaluate the \rhs\ in $O(kn^2)$ flops
as opposed to the $O(kn^3)$ flops that are required if we explicitly
form the product.
If the $A_i$ are bidiagonal
then the costs are $3kn$ flops 
compared with up to $O(kn^2)$ flops if the product is
explicitly formed, since in general the product fills in.

The equality \eqref{Aprod-eq} obviously holds when the $B_i$ are all
nonnegative.
It can also hold because all additions in the
product $A_1A_2\dots A_k$ are of like-signed numbers, so that there is no
cancellation. Important such cases
are when the $A_i$ are nonnegative
and when each $A_i$ has a checkerboard (alternating) sign pattern, 
which can be expressed as
\begin{equation}\label{AiSig}
  A_i = \pm \Sig |A_i| \Sig,  \quad i = 1\colon k,
\end{equation}
where 
\begin{equation}\label{Sig-def}
   \Sig = \diag\bigl(1,-1,1,\dots,(-1)^{n-1}\bigr).
\end{equation}

\begin{theorem}\label{thm.Ainv-prod}
If the matrices 
$A_i$, $i=1\colon k$, satisfy \eqref{AiSig} then
\begin{equation}\label{Aprod-mod2}
 A_1A_2\dots A_k = \pm\Sig |A_1|@|A_2|\dots |A_k| \Sig
\end{equation}
and hence
\begin{equation}\label{Aprod-mod3}
 |A_1A_2\dots A_k| = |A_1|@|A_2|\dots |A_k|.
\end{equation}
\end{theorem}

\begin{proof}
If the \(A_i\) satisfy \eqref{AiSig} then
\begin{equation}\notag
  A_1A_2\dots A_k = \pm\Sig |A_1|\Sig\cdot \Sig |A_2|\Sig \dots \Sig|A_k| \Sig
                  = \pm\Sig |A_1|@|A_2|\dots |A_k| \Sig,
\end{equation}
which is \eqref{Aprod-mod2}, and 
\eqref{Aprod-mod3} follows immediately,
\end{proof}
We conclude that if the $A_i$  are nonnegative or have a \csp\ then we can
compute
$\normo{A_1A_2\dots A_k}$ in $O(kn^2)$ flops.

If $B_1, B_2,\dots,B_k$ are bidiagonal and nonnegative
then from Lemma~\ref{lem.bidiag-inv} it is clear that $B_i^{-1}$ has a
checkerboard sign pattern, that is, it satisfies \eqref{AiSig}. 
Therefore by \eqref{Aprod-mod3},
\begin{equation}\label{Bkinv-prod}
    |B_k^{-1} B_{k-1}^{-1}\dots B_1^{-1}|
    = |B_k^{-1}|@|B_{k-1}^{-1}|\dots |B_1^{-1}|.
\end{equation}

The same is true if the $B_i$ have a \csp.

\begin{theorem}\label{thm.Binv-prod}
Let $B_1,B_2,\dots,B_k\in\R^{\nbyn}$ be nonsingular bidiagonal matrices.
If $B_i$ is nonnegative for all $i$ or has a checkerboard sign pattern for
all $i$ then
\begin{equation}\label{Biprod}
   |B_k^{-1} B_{k-1}^{-1} \dots B_{1}^{-1}| 
  = |B_k^{-1}|@ |B_{k-1}^{-1}| \dots |B_{1}^{-1}| 
  = M(B_k)^{-1} M(B_{k-1})^{-1} \dots M(B_{1})^{-1}. 
\end{equation}
\end{theorem}

\begin{proof}
For nonnegative $B_i$ the result follows from \eqref{Bkinv-prod} on recalling
\eqref{B-MB}.
From \eqref{Binv-ij} it is clear that $B_i$ having a checkerboard sign
pattern is equivalent to either $B_i^{-1}$ or $-B_i^{-1}$
being nonnegative and equal to $M(B_i)^{-1}$,
which gives the second part of the result.
\end{proof}

From \eqref{Biprod} we have 
\begin{equation}\label{Biprod2}
  \normo{ B_k^{-1} B_{k-1}^{-1} \dots B_{1}^{-1} }
  = \normo{ M(B_k)^{-1} M(B_{k-1})^{-1} \dots M(B_{1})^{-1}e },
\end{equation}
and the right-hand side can be computed in $3kn$ flops,
whereas explicitly forming the product on the left 
(using substitutions) costs $3kn^2/2$ flops.
We conclude that when the $B_i$ are nonnegative for all $i$ or all have a
checkerboard sign pattern,
$\kinf(B_1B_2\dots B_k)$ can be computed exactly in $6kn$ flops.
Since $\normi{A} = \normo{A^T}$, the $1$-norm \cn\ can be computed at the
same cost by working with the transpose of the product.

In the case $k = 1$, \eqref{Biprod2} reduces to the result that
$\normo{B^{-1}} = \normo{M(B)^{-1}} = \normo{M(B)^{-1}e}$
\cite[sec.~2]{high86t}.

We can also compute the \cn\ of Skeel \cite{skee79},
\begin{equation}\notag
 \cond(A,x) = \frac{ \normoq{ |A^{-1}|@ |A|@ |x| } }{\normo{x}},
\end{equation}
exactly in $6kn$ flops for $A = B_1B_2\dots B_k$ with nonnegative $B_i$:
\begin{equation}\notag
 \cond(B_1B_2\dots B_k,x) = \frac{\normoq{ M(B_k)^{-1} \dots M(B_1)^{-1} B_1\dots B_k |x|}}
                   {\normo{x}}.
\end{equation}
If the $B_i$ have \csp s then the same formula holds
with $B_1B_2\dots B_k$ replaced by 
     $|B_1| |B_2| \dots |B_k|$.

We will make use of \eqref{Biprod2} for \TN\ matrices in
Section~\ref{sec.TN}.

\section{Linear Systems}\label{sec.linear-systems}

We consider a linear system $A x= b$ in which $A$ is either a product of
bidiagonal matrices or a product of inverses of bidiagonal matrices.
Our interest is in what can be said about the \berr\ and \ferr\ when such a
system is solved in \fpa.

\subsection{Product of Bidiagonal Matrices}

Suppose $A = B_1B_2\dots B_k$ is a product
of $k$ bidiagonal matrices.
We can solve the system  by solving $k$ bidiagonal systems by substitution.
Standard \rea\ \cite[Lem.~8.2]{high:ASNA2} 
shows that the computed $\xhat$ satisfies
\begin{equation}\label{bprod-berrs}
  (B_1 + \dB_1) (B_2 + \dB_2) \dots(B_k + \dB_k)\xhat = b,
  \quad |\dB_i| \le \gu_2 |B_i|, \quad i=1\colon k.
\end{equation}
Hence the residual is
\begin{align*}
  |b - B_1B_2\dots B_k \xhat| 
   &= 
  \bigl|\bigl((B_1+\dB_1)(B_2 + \dB_2) \dots(B_k + \dB_k) - B_1B_2\dots B_k\bigr)
  \xhat@\bigr|\\
  &\le \bigl( (1+\gu_2)^k-1 \bigr) |B_1|@|B_2|\dots |B_k|@ |\xhat|,
\end{align*}
by Lemma~\ref{lem.matprod-comp}.
If the $B_i$ are all nonnegative or,
by Theorem~\ref{thm.Ainv-prod}, if they 
have a checkerboard sign pattern,
then the bound becomes 
\begin{equation}\label{res0}
   |b - A\xhat| \le \bigl( (1+\gu_2)^k-1 \bigr) |A|@ |\xhat|
                 = \bigl( 2ku + O(u^2) \bigr) |A|@ |\xhat|,
\end{equation}
which shows that the \crbe\ is small---an ideal \berr\ result.
We note that this result has 
used the triangularity of the $B_i$ 
but not their bidiagonal structure 
(except through the constant in \eqref{bprod-berrs}).

To obtain a \ferr\ bound we rewrite \eqref{bprod-berrs} as
\begin{equation*}
   \xhat = 
  (B_k+\dB_k)^{-1}(B_{k-1} + \dB_{k-1})^{-1} \dots(B_1 + \dB_1)^{-1} b.
\end{equation*}
Then 
\begin{align}\notag
 |\xhat - x| &\le \bigl|(B_k+\dB_k)^{-1}(B_{k-1} + \dB_{k-1})^{-1} \dots(B_1 +
               \dB_1)^{-1} - B_k^{-1}  B_{k-1}^{-1} \dots B_1^{-1}\bigr| |b|\\
            &\le \bigl((1+\tau)^k-1\bigr) |B_k^{-1}| @|B_{k-1}^{-1}|@
             \dots |B_1^{-1}| @ |b|
\end{align}
by Theorem~\ref{thm.bidiag-prod-inv-pert}, where
\begin{equation}\label{tau2}
   \tau = \frac{(2n-1)\gu_2}{1 - (2n-1)\gu_2}.
\end{equation}
If the $B_i$ are all nonnegative or have a checkerboard sign pattern
then by Theorem~\ref{thm.Binv-prod}
this inequality becomes
\begin{equation}\label{ferrb}
 |\xhat - x| \le \bigl(2k(2n-1) u + O(u^2)\bigr) |A^{-1}| @|b|.
\end{equation}
The bound \eqref{ferrb}
is a strong \ferr\ bound because it is the same as a bound for
the change in $x$ induced by a small \cmp\ relative \pert\ of 
of $b$: $b\to b+\db$ with $|\db| \le 4knu |b|$~\cite[Thm.~7.4]{high:ASNA2}.

\subsection{Product of Inverses of Bidiagonal Matrices}

Now suppose that it is $A^{-1}$
rather than $A$ that is a product of bidiagonal matrices:
$A^{-1} = B_1B_2\dots B_k$.
Now we solve $Ax = b$ by forming 
$x = A^{-1}b = B_1B_2\dots B_kb$ and the computed $\xhat$ satisfies
\begin{equation}\label{xhat1}
    \xhat = (B_1+\dB_1)(B_2 + \dB_2) \dots(B_k + \dB_k) b,
    \quad |\dB_i| \le \gu_2 |B_i|, \quad i=1\colon k.
\end{equation}
Then the \ferr\ is 
\begin{align}
 |\xhat - x| &= \bigl|\bigl((B_1+\dB_1)(B_2 + \dB_2) \dots(B_k + \dB_k) - B_1B_2\dots
                 B_k)\bigr) b\bigr|, \notag \\
         &\le \bigl( (1+\gu_2)^k - 1\bigr) |B_1|@ |B_2| \dots |B_k|@|b|,
         \label{xhatmx2}
\end{align}
by Lemma~\ref{lem.matprod-comp}.
If the $B_i$ are all nonnegative or 
have a checkerboard sign pattern then
by Theorem~\ref{thm.Ainv-prod},
$ |B_1|@ |B_2| \dots |B_k| = |B_1 B_2 \dots B_k|$, so 
\begin{equation}\label{ferr-b}
 |\xhat - x| \le \bigl( (1+\gu_2)^k - 1\bigr) |A^{-1}|@|b|.
\end{equation}

Now we turn to the residual.  Note first that by \eqref{xhat1},
\begin{equation}\notag
 b = (B_k + \dB_k)^{-1} (B_{k-1} + \dB_{k-1})^{-1} \dots(B_1+\dB_1)^{-1} \xhat.
\end{equation}
Hence
\begin{equation}\notag
 |b-A\xhat| = \bigl| \bigl[ (B_k + \dB_k)^{-1} (B_{k-1} + \dB_{k-1})^{-1}
                 \dots(B_1+\dB_1)^{-1} - B_k^{-1} B_{k-1}^{-1}\dots
                 B_1^{-1} \bigr] \xhat\bigr|
\end{equation}
and by Lemma~\ref{lem.matprod-comp} and 
Theorem~\ref{thm.bidiag-prod-inv-pert} we obtain, 
with $\tau$ given by \eqref{tau2},
\begin{align*}
  |b - A\xhat| &\le \bigl( ( 1+ \tau )^k -1 \bigr) 
                |B_k^{-1}|@ |B_{k-1}^{-1}| \dots |B_1^{-1}|@ |\xhat|\\
       & = \bigl(2k(2n-1)u + O(u^2)\bigr) |B_k^{-1}|@ |B_{k-1}^{-1}| \dots |B_1^{-1}|@ |\xhat|.
\end{align*}
If the $B_i$ are all nonnegative or have a checkerboard sign pattern 
then by Theorem~\ref{thm.Binv-prod} this bound can be written
\begin{equation}\label{res1}
 |b - A\xhat| \le \bigl(2k(2n-1)u + O(u^2)\bigr) |A|@ |\xhat|,
\end{equation}
which again shows a small \crbe.

Our conclusion is that  whether it is $A$ or $A^{-1}$ that is a product of
bidiagonal matrices we have the same satisfactory form of \ferr\
bounds \eqref{ferrb} and \eqref{ferr-b} 
and residual bounds \eqref{res0} and \eqref{res1} when 
the $B_i$ are all nonnegative or have a \csp.

\subsection{Application to Vandermonde Systems}\label{sec.vand}

An application of these results is to the Bj\"orck--Pereyra \alg\
for solving a Vandermonde system $Vy = b$ in $O(n^2)$ flops \cite{bjpe70},
where $V = (x_j^{i-1}) \in \C^{\nbyn}$ for given points $x_i\in\C$.
This \alg\ uses a \fact\ of $V^{-1}$ into a product of $2n-2$ bidiagonal
matrices $B_{2n-2}, \dots, B_1$ given in terms of the points $x_i$.
When $0\le x_1 < x_2 < \cdots < x_n$ the bidiagonal factors have positive
diagonal and nonpositive off-diagonal elements.
Therefore the $B_i$ have a checkerboard sign pattern and so 
$|B_{2n-2}| \dots |B_1| = |B_{2n-2} \dots B_1| = |A^{-1}|$ by 
\eqref{Bkinv-prod}.
From \eqref{ferr-b} and \eqref{res1} we have 
\begin{align*}\notag
 |\yhat-y| &\le \bigl(2(2n-2)u + O(u^2) \bigr) |V^{-1}| @ |b|,\\
 |b - V\yhat| &\le \bigl(2(2n-2)(2n-1)u + O(u^2)\bigr) |V|@ |\yhat|,
\end{align*}
which reproduce \cite[Thm.~2.3]{high87v}
and the monomial case of 
\cite[Cor.~4.1]{high90v}, \resp.
Since $V^{-1}$ has a checkerboard sign pattern, if $(-1)^ib_i\ge0$
then $|V^{-1}| @ |b| = |V^{-1}b| = |y|$,
and $\yhat$ therefore has a small \cmp\ relative error.
The analysis in \cite{high90v} makes use of the bidiagonal \fact,
but that in \cite{high87v} does not.

\subsection{Application to Pascal Systems}\label{sec.pascal0}

We give a numerical illustration of the use of the bidiagonal \fact\ for  
solving the linear system $P_nx = b$, where $P_n$ is the \spd\ $\nbyn$ 
Pascal matrix with 
\begin{equation}\label{pascal-mat}
   p_{ij} = {i+j-2 \choose j-1} = \frac{ (i+j-2)! }{ (i-1)! (j-1)! }
\end{equation}
and $b = e_n/n$, where $e_n$ is the $n$th unit vector.
The Pascal matrix 
has a known \fact\ as a product of $2n-1$ 
bidiagonal matrices, as we explain in section~\ref{sec.pascal}.
We solve the system using the bidiagonal \fact, solving the bidiagonal 
systems by substitution.
We also solve the system for the explicitly formed $P$ using the
MATLAB backslash operator
(which exploits the symmetric positive definiteness of $P_n$ 
but not its bidiagonal \fact).
The working precision is double precision, 
with $u \approx 1.1 \times 10^{-16}$.
Table~\ref{table.pascal-err} shows the relative errors 
$\normo{x-\xhat}/\normo{x}$,
for which we take as the exact solution $x$
the solution computed at a precision of 500 decimal digits using the
Multiprecision Computing Toolbox \cite{adva-mct} and then rounded to double
precision.
We restrict to $n \le 25$ to ensure that $P$ is 
exactly representable at the working precision.
We see that substitution with the bidiagonal \fact\
yields errors of $O(u)$ that satisfy the bound
\eqref{ferrb}, whereas
the MATLAB backslash function produces much larger errors, which usually exceed
\eqref{ferrb}.

\begin{table}
\caption{Relative errors for the computed solution to 
         a linear system $P_n x = b$ with $P_n$ the $\nbyn$ Pascal matrix.}
\label{table.pascal-err}

\begin{center}
\begin{tabular}{cccc}\toprule
          & \multicolumn{2}{c}{Relative errors}\\\cmidrule(lr){2-3} 
$n$ &Bidiagonal \fact &  \verb"P\b" & Error bound \eqref{ferrb} \\\midrule
  5 &9.25e-17         &  9.25e-16   & 7.99e-15                  \\
 10 &1.50e-16         &  4.94e-9    & 3.80e-14                  \\
 15 &6.36e-17         &  1.05e-3    & 9.02e-14                  \\
 20 &1.34e-16         &  3.12e-12   & 1.65e-13                  \\
 25 &1.68e-16         &  2.76e-11   & 2.61e-13                  \\
\bottomrule
\end{tabular}
\end{center}

\end{table}

\section{Totally Nonnegative Matrices}\label{sec.TN}

A matrix $A\in\mathbb{R}^{n\times n}$ is 
totally nonnegative if every
minor (determinant of a square submatrix) is nonnegative
and totally positive if every minor is positive.
We will need the following key result, which is a
direct consequence of the Binet--Cauchy theorem on
determinants \cite[sec.~0.8.7]{hojo13}.

\begin{theorem}\label{thm.TN-prod}
If $A,B\in\mathbb{R}^{n\times n}$ are totally nonnegative 
then so is $AB$.
\end{theorem}

Bidiagonal matrices play a key role in the theory of \TN\ matrices.
Indeed a nonnegative bidiagonal matrix is \TN.
In the proof of this result we will need 
the elementary lower bidiagonal matrix
\begin{equation}\label{Lk-elt}
 L_k(\ell_{k+1,k}) = I + \ell_{k+1,k} e_{k+1}e_k^T,
\end{equation}
which differs from the identity matrix only in the $(k+1,k)$ position,
which contains $\ell_{k+1,k}$.

\begin{theorem}\label{thm.tn-bidiag}
A bidiagonal matrix $B\in\R^{\nbyn}$ 
with nonnegative elements is \TN.
\end{theorem}

\begin{proof}
\Wlogy\ we take $B$ to be lower bidiagonal.
We first assume that $B$ is nonsingular.
Since $0\ne \det(B) = b_{11}b_{22}\dots b_{nn}$,
the $b_{ii}$ are all positive, so with $D = \diag(b_{ii})$ 
and $\ell_{i+1,i} = b_{i+1,i} / b_{i+1,i+1} \ge 0$, $i=1\colon n-1$,
we can write 
\begin{equation}\label{B-DL}
 B = D 
     \begin{bmatrix} 
        1         &           &        &              & \\
        \ell_{21} & 1         &        &              & \\
                  & \ell_{32} & \ddots &              & \\
                  &           & \ddots & \ddots       & \\ 
                  &           &        & \ell_{n,n-1} & 1 
     \end{bmatrix}
     \equiv DL.
\end{equation}
Since $D$ is clearly \TN, 
by Theorem~\ref{thm.TN-prod} it suffices to show that $L$ is \TN.

For $n = 4$ we have 
\begin{equation*}
 L = \begin{bmatrix} 
        1         &           &           &   \\
        \ell_{21} & 1         &           &   \\
                  & \ell_{32} & 1         &   \\
                  &           & \ell_{43} & 1 \\
     \end{bmatrix}
 = 
     \begin{bmatrix} 
        1         &           &           &   \\
        \ell_{21} & 1         &           &   \\
                  &           & 1         &   \\
                  &           &           & 1 \\
     \end{bmatrix}
     \begin{bmatrix} 
        1         &           &           &   \\
                  & 1         &           &   \\
                  & \ell_{32} & 1         &   \\
                  &           &           & 1 \\
     \end{bmatrix}
     \begin{bmatrix} 
        1         &           &           &   \\
                  & 1         &           &   \\
                  &           & 1         &   \\
                  &           & \ell_{43} & 1 \\
     \end{bmatrix},
\end{equation*}
and this \fact\ clearly generalizes to 
\begin{equation}\label{nonnegL-fact}
L = L_1(\ell_{21})
     L_2(\ell_{32}) \dots
     L_{n-1}(\ell_{n,n-1}),
\end{equation}
where $L_k(\ell_{k+1,k})$ is the elementary lower bidiagonal
matrix \eqref{Lk-elt}.
It is easy to see that $L_k(\ell_{k+1,k})$ is \TN\ for all $k$, 
so $L$ is \TN\ by Theorem~\ref{thm.TN-prod}.

If $B$ is singular then consider the bidiagonal matrix
$B(\eps) = B + \eps I$, which is nonsingular for $\eps > 0$.
By the argument above, $B(\eps)$ is \TN\ for $\eps > 0$.
Any minor of $B(\eps)$ is the determinant of a submatrix of $B(\eps)$,
which is a \py\ in $\eps$, so it is continuous in $\eps$.
This minor is nonnegative for all $\eps > 0$
and so must remain nonnegative in the limit as $\eps \to 0$.
Therefore $B = B(0)$ is \TN.
\end{proof}

Even if $B$ is not \TN, there is a an associated \TN\ matrix.

\begin{theorem}\label{thm.MTinv-TN}
If $B\in\R^{\nbyn}$ is nonsingular and bidiagonal then $M(B)^{-1}$ is \TN.
\end{theorem}

\begin{proof}
Assuming that $B = L$ is lower bidiagonal,
by \eqref{B-DL} and \eqref{nonnegL-fact},
\begin{equation}\notag
M(B) = M(DL) = |D| M(L) 
        = |D| L_1(-|\ell_{21}|)
        L_2(-|\ell_{32}|) \dots
        L_{n-1}(-|\ell_{n,n-1}|)
\end{equation}
and $L_k(-|\ell_{k+1,k}|)^{-1} = L_k(|\ell_{k+1,k}|)$,
so 
$M(B)^{-1} =  L_{n-1}(|\ell_{n,n-1}|)
              L_{n-2}(|\ell_{n-1,n-2|})
              \dots L_1(|\ell_{21}|) \* |D|^{-1}$,
which is a product of \TN\ matrices and hence is \TN.
\end{proof}

The next result
shows that any nonsingular \TN\ matrix can be written as a product of 
nonnegative bidiagonal matrices.

\begin{theorem}\label{thm.TN-ALUfact}
A nonsingular matrix 
$A\in\mathbb{R}^{n\times n}$ is totally nonnegative if and only if it 
can be factorized as
\begin{equation}\label{TN-ALUfact}
   A = L_{n-1} L_{n-2} \dots L_1 D U_1 U_2 \dots U_{n-1}, 
\end{equation}
where $D$ is a diagonal matrix with positive diagonal entries
and $L_i$ and $U_i$ are unit lower and unit upper bidiagonal matrices,
respectively, with the first $i-1$ entries along the subdiagonal 
of $L_i$ and $U_i^T$ zero and the rest nonnegative.
\end{theorem}

The \fact\ \eqref{TN-ALUfact} is essentially an \LUf\ in which $L$ and $U$
have been factorized into a product of specially structured nonnegative
bidiagonal matrices.

Theorem~\ref{thm.TN-ALUfact} is from 
Gasca and Pe\~na \cite[Thm.~4.2]{gape96}.
Fallat and Johnson \cite[sec.~2.0]{fajo11} summarize
the history of different forms of this \fact.

Since the bidiagonal matrices in the \fact\ \eqref{TN-ALUfact} are all
nonnegative, by \eqref{Biprod2} we have 
\begin{equation}\label{Biprod3}
 \normo{A^{-1}}  
  = \normo{ M(U_{n-1})^{-1} \dots\ M(U_1)^{-1} D^{-1}
            M(L_1)^{-1}\dots\ M(L_{n-1})^{-1}e },
\end{equation}
and so we can compute $\normo{A^{-1}}$ by $2(n-1)$ substitutions
in $O(n^2)$ flops for any nonsingular 
\TN\ matrix given the \fact\ \eqref{TN-ALUfact}.

Let $\chat = \fl(\normo{A^{-1}})$.
Taking $\infty$-norms in \eqref{ferrb} with $b = e$ gives,
using the triangle inequality,
\begin{equation}\label{conderr}
 \frac{|\chat - \normo{A^{-1}}|}{\normo{A^{-1}}} \le dn^2u
\end{equation}
for a modest constant $d$.
Therefore $\chat$ is highly accurate,
essentially because there is no cancellation in evaluating \eqref{Biprod3}:
all additions are of nonnegative quantities.
Standard methods for evaluating $\normo{A^{-1}}$ for general $A$ only satisfy 
$|\chat - \normo{A^{-1}}|/\normo{A^{-1}} \le cn^3 \kinf(A)u$, 
which is the best that can be expected in general because the \cn\ of 
$\kinf(A)$ is $\kinf(A)$~\cite{high95cn}.

To obtain $\kinf(A)$ we need $\normo{A}$, which can either be computed from 
$A$ if it is explicitly known, or from 
$\normo{A} = \normo{L_{n-1} L_{n-2} \dots L_1 D U_1 U_2 \dots U_{n-1} e}$ 
otherwise. We summarize the computations in an \alg.

\begin{algorithm}\label{alg.fastTNcond}
This \alg\ computes $c = \kinf(A)$ for a \TN\ matrix $A$ given the \fact\
\eqref{TN-ALUfact}.
\begin{code}
If \= $A$ is explicitly known\\ 
  \> $\a = \normo{A}$\\
else\\
  \> $\a = \normo{L_{n-1} L_{n-2} \dots L_1 D U_1 U_2 \dots U_{n-1} e}$\\
end\\
Compute
$\b  = \normo{ M(U_{n-1})^{-1} \dots\ M(U_1)^{-1} D^{-1}
            M(L_1)^{-1}\dots\ M(L_{n-1})^{-1}e }$\nonumberbreak
by substitutions.\\
$c = \a \b$
\end{code}
\end{algorithm}

How do we obtain the parameters in the \fact\ \eqref{TN-ALUfact}?
In some cases they are known from the construction of the matrix.
Formulas are known for \TP\ Vandermonde matrices and Cauchy matrices
\cite[eqs.~(3.5), (3.6)]{koev05}
and a variety of Vandermonde-type matrices \cite{dkmm23}.
For \TP\ matrices determinantal formulas for the parameters are available
\cite[Prop.~3.1]{koev05}.
Assuming the determinants can be computed accurately,
in all these cases the parameters can be evaluated to high relative
accuracy.
and so in view of Theorem~\ref{thm.bidiag-prod-inv-pert} the errors in the
evaluation of the parameters do not affect the form of the bound
\eqref{conderr}.

We give two numerical experiments in MATLAB to illustrate the accuracy of
the \cn\ evaluation.  
We take as the exact \cn\ the 
one computed at a precision of 500 decimal digits using the
Multiprecision Computing Toolbox
\cite{adva-mct} and then rounded to double precision.

First, in Table~\ref{table.Hilbert} we show the relative errors in
computing the $\infty$-norm \cn\ of the Hilbert matrix $H_n$, 
which has $(i,j)$ element $1/(i+j-1)$ and is \TP.
The parameters in the bidiagonal \fact\ \eqref{TN-ALUfact} are computed
using the function \verb"TNCauchyBD" from 
the TNTool 
toolbox.\footnote{\url{http://www.math.sjsu.edu/~koev/software/TNTool.html}}
We see that even extremely large \cn s are obtained to high accuracy.

\begin{table}
\caption{Condition numbers and relative errors for the Hilbert matrix.}
\label{table.Hilbert}

\begin{center}
\begin{tabular}{ccc}\toprule
$n$       & $\kinf(H_n)$& Relative error for \Alg~\ref{alg.fastTNcond}\\ \midrule
  4       & 2.84e4    & 1.28e-16  \\
  8       & 3.39e10   & 2.25e-16  \\
 16       & 5.06e22   & 3.67e-17  \\
 32       & 1.36e47   & 1.75e-15  \\
 64       & 1.10e96   & 1.77e-15  \\
\bottomrule
\end{tabular}
\end{center}

\end{table}

Next we consider the Pascal matrix \eqref{pascal-mat}, which is \TP\
\cite[Ex.~0.1.6]{fajo11}.
Since this matrix is exactly representable at the working precision for 
$n$ up to around $25$, we can 
compare \Alg~\ref{alg.fastTNcond}
with the MATLAB \verb"cond" function.
We see from the results in Table~\ref{table.Pascal}
that the MATLAB function loses accuracy as $n$ increases while
\Alg~\ref{alg.fastTNcond} returns a result correct to the working precision.

\begin{table}
\caption{Condition numbers and relative errors for the Pascal matrix.}
\label{table.Pascal}

\begin{center}
\begin{tabular}{cccc}\toprule
    &            & \multicolumn{2}{c}{Relative errors}            \\\cmidrule(lr){3-4} 
$n$ & $\kinf(P_n)$  & \Alg~\ref{alg.fastTNcond} & \verb"cond(P_n,inf)"\\ \midrule
  5 & 1.56e4      & 0.00                      & 0.00              \\
 10 & 8.13e9      & 0.00                      & 1.49e-11          \\
 15 & 5.77e15     & 0.00                      & 2.19e-8           \\
 20 & 4.50e21     & 4.66e-17                  & 3.41e-4           \\
 25 & 3.81e27     & 1.70e-17                  & 3.17e-2           \\
\bottomrule
\end{tabular}
\end{center}

\end{table}

Another use of the \fact\ 
of Theorem~\ref{thm.TN-ALUfact}
is to construct \TN\ matrices by choosing the $n^2$ parameters 
that make up the $L_i$, $D$, and the $U_i$.
The function call
\begin{lstlisting}
A = anymatrix('core/totally_nonneg',X)
\end{lstlisting}
in the Anymatrix toolbox \cite{himi21} constructs an
$\nbyn$ \TN\ matrix $A$ from parameters given in the $\nbyn$ matrix
$X$, whose format is as suggested in \cite[sec.~4]{koev05}.  
The Pascal matrix is generated when \lstinline"X = ones(n)".
In a call
\begin{lstlisting}
A = anymatrix('core/totally_nonneg',n)
\end{lstlisting}
the parameters are chosen randomly,
and this is a convenient way to generate random \TN\ matrices.

Koev \cite[sec.~7]{koev05}, \cite{koev07}
shows that small relative
changes in the parameters in the \fact\ \eqref{TN-ALUfact}
produce small relative changes in the 
the determinant, the \evals, and the \svals.
In \cite{koev05} he develops \alg s for accurate computation of
\evals\ and the SVD of nonsingular \TN\ matrices,
given an accurate bidiagonal \fact,
by carrying out transformations on the bidiagonal
\fact\ in such a way that no subtractions occur.

For later use, we note a useful theorem about the \evals\ of a
\TN\ matrix \cite[Thm.~3.3]{fgj00}.

\begin{theorem}\label{thm.TN-eig-irred-sing}
If $A\in\mathbb{R}^{n\times n}$ is totally nonnegative
and irreducible then  its \evals\ are real and nonnegative 
and the positive \evals\ are distinct.
\end{theorem}

Note that the irreducibility requirement in the theorem means that it
cannot be applied to \tri\ matrices, so there is no contradiction to the
fact that the \TN\ matrix $\smatrix{1 & 1 \\ 0 & 1}$ (for example) has
repeated nonzero \evals.

\section{Matrix Functions and Polynomial Evaluation and Interpolation}\label{sec.matrix-functions}

Bidiagonal matrices are intimately connected with \py\ evaluation and
interpolation.
Horner's method for evaluating a \py\ at a point $\a$ 
can be expressed as the solution 
of a linear system with coefficient matrix $T_n(-\a)$
\cite[sec.~5.2]{high:ASNA2}, where $T_n$ is defined in \eqref{Tnth}.
Premultiplying a vector by $T_n(-1)^T$ corresponds to forming a backward 
difference, and a subsequent \mult\ by a diagonal matrix 
yields divided differences \cite[sec.~5.3]{high:ASNA2}.
In fact, an explicit formula for a function of a bidiagonal matrix is
available in terms of divided differences.
Recall that divided differences of a function $f$ 
at points $x_k$ are defined recursively by
(see, e.g. \cite[Chap.~2]{code80} or \cite[sec.~B.16]{high:FM} )
\begin{align}
      f[x_k] &= f(x_k), \nonumber\\[\medskipamount]
  f[x_0,x_1,\dots,x_{k+1}] &=
  \begin{cases}
  \dfrac{ f[x_1,x_2,\dots,x_{k+1}]  - f[x_0,x_1,\dots,x_k] }%
       { x_{k+1} - x_0 }, & x_0\ne x_{k+1},\\
         \dfrac{ \strut f^{(k+1)}(x_{k+1}) }{\strut (k+1)!},
         & x_0 = x_{k+1}, \qquad \end{cases}   
    \label{dd-conf}
\end{align}
where, since $f[x_1,x_2,\dots,x_{k+1}]$ does not depend on the order of its
arguments, we assume without loss of geniality that equal points are
contiguous.

\begin{theorem}\label{thm.fB}
If $B\in\C^{\nbyn}$ is upper bidiagonal then $F = f(B)$ is upper \tri\ with
$f_{ii} = f(t_{ii})$ and
\begin{equation}\label{fB-eq}
      f_{ij} = b_{i,i+1} b_{i+1,i+2} \dots b_{j-1,j}\,
               f[b_{ii},b_{i+1,i+1}, \dots, b_{jj}], \quad j > i.
\end{equation}
\end{theorem}

\begin{proof}
The formula \eqref{fB-eq}
is a special case of the formula for $f(T)$,
where $T$ is upper \tri,
given in Davis \cite{davi73}, Descloux \cite{desc63a},
and Van Loan \cite{vanl75}.
\end{proof}

Lemma~\ref{lem.bidiag-inv} is the special case of Theorem~\ref{thm.fB}
with $f(x) = 1/x$.
Since $f[\l,\l,\dots,\l] = f^{(n-1)}(\l)/(n-1)!$, 
another special case is the formula for a 
function of an $\mbym$ Jordan block~\cite[sec.~1.2]{high:FM}
\begin{equation}\label{f-Jblock}
 f\left(
     \Bmatrix{ \lambda   & 1         &          &           \cr
                         & \lambda   & \ddots   &           \cr
                         &           & \ddots   &    1      \cr
                         &           &          & \lambda   \cr } \right)
         =    \Bmatrix{ f(\l) & f'(\l) &  \dots  & \DS\frac{f^{(m-1)}(\l)}{(m-1)!} \cr
                     & f(\l)  & \ddots  & \vdots \cr
                     &          & \ddots  & f'(\l) \cr
                     &          &         & f(\l) \cr}.
\end{equation}
Yet another special case is
\begin{equation}\notag
  f\left(\Bmatrix{ \l_1 & 1    &        &      \cr
                        & \l_2 & \ddots &      \cr
                        &      & \ddots & 1    \cr
                        &      &        & \l_n \cr} \right)_{1n}
    = f[\l_1,\l_2,\dots,\l_n],
\end{equation}
which is a result of Opitz~\cite{opit64}
and is used in computing divided differences of the exponential
by McCurdy, Ng, and Parlett \cite{mnp84}.

A natural question is whether a function of a nonnegative bidiagonal
matrix is \TN.
For the exponential, the answer is yes.

\begin{theorem}\label{thm.exp-TNB}
If $B\in\R^{\nbyn}$ is a nonnegative bidiagonal matrix then 
$\eu^B$ is \TN.
\end{theorem}

\begin{proof}
Consider the formula \cite[sec.~10.1]{high:FM}
$\eu^A = \lim_{m\to \infty} (I + A/m)^m$, valid for any $A$,
where $m\in\Z$.
For nonnegative bidiagonal $B$,
$I + B/m \ge 0$ for all $m > 0$, 
so by Theorem~\ref{thm.tn-bidiag} $I + B/m$ is \TN\
and therefore $X_m = (I + B/m)^m$ is \TN\ for all $m>0$
by Theorem~\ref{thm.TN-prod}.
Suppose that 
$\lim_{m\to\infty}X_m$ 
is not \TN, so that 
some submatrix with indices $(\a,\b)$ has negative determinant.
Let $x_m = \det(X_m(\a,\b))$.
Then $\lim_{m\to\infty} x_m < 0$ but $x_m>0$ for all $m$,
which is a contradiction,
so $\eu^B$ is \TN.
\end{proof}

Note that Theorem~\ref{thm.exp-TNB} does not generalize to wider
bandwidths, as the example
\begin{equation}\notag
      \exp\left( 
       \Bmatrix{1 & 1 & 1 \cr
                  & 1 & 1 \cr
                  &   & 1 }
       \right)
    =  \Bmatrix{\eu & \eu & 3\eu/2 \cr
                  & \eu & \eu \cr
                  &   & \eu }
\end{equation}
shows, since the $(1\colon 2, 3\colon 4)$  submatrix has negative determinant.

\section{Upper Triangular Toeplitz matrices}\label{sec.toep-upper}

Upper \tri\ Toeplitz matrices $T\in\C^{\nbyn}$ can be written in the form
\begin{equation}\notag
T =  \Bmatrix{ t_0 & t_1 & \dots  & t_{n-1}\cr
                   & t_0 & \ddots & \vdots \cr
                   &     & \ddots & t_1   \cr
                   &     &        & t_0\cr} 
                 = \sum_{j=1}^n t_{j-1} N^{j-1},
\end{equation}
where $N$ is upper bidiagonal with a superdiagonal of ones:
\begin{equation}\notag
     N = \Bmatrix{0 & 1      &        &   \cr
                    & 0      & \ddots &   \cr
                    &        & \ddots &   1  \cr
                    &        &        &   0}, 
\end{equation}
Note that $N^n = 0$.
It follows that 
the product of two upper \tri\ Toeplitz matrices is again 
upper \tri\ Toeplitz
and 
that 
upper \tri\ Toeplitz matrices commute.
Furthermore, since 
$f(T)$ is a \py\ in $T$, it follows that 
$f(T)$ is also upper \tri\ and Toeplitz.
Note that as a special case,
if $B$ is a Toeplitz bidiagonal matrix with 
$b_{ii} = b$ and $b_{i,i+1} = c$ then
Theorem~\ref{thm.fB} gives
$f(B)_{ij} = c^{j-i} f[b,b,\dots,b]
           = c^{j-i} f^{(j-i)}(b)/(j-i)!$,
of which \eqref{f-Jblock} is a special~case.

\section{Exploiting Factorizations Into Products of Bidiagonal Matrices}\label{sec.expl-fact}

In this section we show how \fact s involving
bidiagonal matrices or their inverses can provide valuable information about
particular matrices.

\subsection{The Frank Matrix}\label{sec.frank}

In 1958 Frank~\cite{fran58} reported that his \alg s had difficulty
computing accurately the smaller \evals\ of the 
$\nbyn$ upper Hessenberg matrix 
\begin{equation}\notag
  F_n = 
  \left[\begin{array}{*{6}{c}}
  n      & n-1    & n-2    & \dots  & 2      & 1 \\
  n-1    & n-1    & n-2    & \dots  & 2      & 1 \\
  0      & n-2    & n-2    & \dots  & 2      & 1 \\[-4pt]
  \vdots & 0      & \ddots & \ddots & \vdots & 1 \\[-5pt]
  \vdots & \vdots & \dots  &   2    & 2      & 1 \\
  0      & 0      & \dots  &   0    & 1      & 1 \\
  \end{array}\right].
\end{equation}
Wilkinson~\cite[sec.~8]{wilk60}
\cite[pp.~92--93]{wilk65} showed that the difficulties
are caused by
the sensitivity of the eigenvalues to perturbations in the matrix, 
which can be measured by the \cn\ of a\ simple \eval\ $\l$:
$\k_2(\l) = \normt{y}\normt{x}/|y^*x|$, where $x$ and $y$ are right and
left \evecs, \resp, corresponding to $\l$.
The \evals\ are known to be real and positive,
and they can be expressed in terms of the zeros of Hermite polynomials
\cite{eber71a}, \cite{vara86}.
However, in none of these references is it shown that the \evals\ are
distinct, which is necessary for the \eval\ \cn s to be defined.

If we subtract row $k+1$ from row $k$ for $k =1\colon n-1$, we obtain a
lower bidiagonal matrix.  For $n=4$ this transformation can be written
\begin{equation}\notag
  \Bmatrix{ 1 & -1 &   &  \cr
              &  1 &-1 &  \cr
              &    & 1 &-1 \cr
              &    &   & 1 \cr}
  \Bmatrix{ 4 &  3 & 2 & 1\cr
            3 &  3 & 2 & 1\cr
              &  2 & 2 & 1\cr
              &    &  1& 1 \cr}
= 
  \Bmatrix{ 1 &    &   &  \cr
            3 &  1 &   &  \cr
              &  2 & 1 &  \cr
              &    & 1 & 1\cr},
\end{equation}
and in general we have 
\begin{align}
      F_n = 
            T_n(-1)^{-1}
            \Bmatrix{  1 &    &        &        &   \cr
                     n-1 &  1 &        &        &   \cr
                         & n-2&  1     &        &   \cr
                         &    &  \ddots& \ddots &       \cr
                         &    & \phantom{n-3}   &   1    & 1     \cr}
          \equiv T_n(-1)^{-1}L,
         \label{frank-fact}
\end{align}
where $T_n$ is defined in \eqref{Tnth}.
This is equivalent to a \fact\ noted by Rutishauser~\cite[sec.~9]{ruti68}.
Note that this is a $UL$ \fact, not an \LUf,
and it takes advantage of the rank-$1$ nature of the upper triangle of $F_n$.
This \fact\ shows that the inverse
$F_n^{-1} = L^{-1}T_n(-1)$ is lower Hessenberg with
integer entries
and that $\det(F_n) = 1$.
Furthermore,
$L$ is \TN\ by Theorem~\ref{thm.tn-bidiag}
and $T_n(-1)^{-1} = M(T_n(-1))^{-1}$ is \TN\ by Theorem~\ref{thm.MTinv-TN}, 
so $F_n$ is the product of two \TN\ matrices and so is \TN\ 
by Theorem~\ref{thm.TN-prod}---a property that to our knowledge 
has not previously been noted.
Since $F_n$ is nonsingular, irreducible
(being upper Hessenberg with nonzero subdiagonal), and \TN\ it follows from 
by Theorem~\ref{thm.TN-eig-irred-sing} that $F_n$ has distinct \evals.
The distinctness of the \evals\ also follows from some rather lengthy
analysis of the characteristic \py\ in \cite[Thm.~2.5]{meba21}.

Frank discussed two matrices in his paper.
The other matrix is obtained from $A_n = (\min(i,j)) \in\R^{\nbyn}$ 
by taking the rows and columns in reverse order.  
We will focus on $A_n$.
For example,
\begin{equation}\notag
    A_4 = \Bmatrix{ 
         1 & 1 & 1 & 1\cr 1 & 2 & 2 & 2\cr 1 & 2 & 3 & 3\cr 1 & 2 & 3 & 4 }.
\end{equation}
The determinant, the inverse, and the \evals\ of $A_n$ can all be easily 
found by constructing a \fact\ involving a bidiagonal matrix.
Consider subtracting row $k-1$ from row $k$  for $k = n:-1:2$.
For $A_4$ this yields
\begin{equation}\notag
     \Bmatrix{ 
         1  &    &   &  \cr 
         -1 & 1 &   &  \cr 
            &-1 & 1 &  \cr 
            &   & -1& 1 } 
    \Bmatrix{ 
         1 & 1 & 1 & 1\cr 1 & 2 & 2 & 2\cr 1 & 2 & 3 & 3\cr 1 & 2 & 3 & 4 }
        = \Bmatrix{
         1 & 1 & 1 & 1\cr   & 1 & 1 & 1\cr   &   & 1 & 1\cr   &   &   & 1 }.
\end{equation}
In general, $T_n(-1)^T A_n = U$,
where $U$ is the upper \tri\ matrix of $1$s.
Hence $A_n = T_n(-1)^{-T} U$, which is a Cholesky \fact\
$A_n = U^TU$ since $T_n(-1)^{-1} = U$,
which shows that $A_n$ is \spd.
Furthermore, $\det(A) = \det(U)^2 = 1$ and 
$A_n^{-1} = U^{-1}U^{-T} = T_n(-1) T_n(-1)^T$, 
which is tridiagonal since $T_n$ is upper bidiagonal.
Now 
$T_n(-1)^{-1}$ is \TN, as noted above; hence 
$A_n$ is the product of two \TN\ matrices and therefore is \TN.
By Theorem~\ref{thm.TN-eig-irred-sing}, the \evals\ of $A_n$ are distinct.
In fact, $A_n^{-1}$ is the almost-Toeplitz tridiagonal matrix
\begin{equation}\notag
 A_n^{-1} = 
      \Bmatrix{
                       2   & -1&        &        &   \cr
                       -1  & 2 & -1     &        &   \cr
                           & -1& \ddots & \ddots &   \cr
                           &   & \ddots & 2      & -1\cr
                           &   &        & -1     & 1},
\end{equation}
and its \evals\ are 
\cite{elli53}, \cite[Chap.~7]{grka69}
(and as given by Frank)
\begin{equation}\notag
  \mu_k = 2 \biggl (1 + \cos\biggl( \frac{2k \pi}{2n+1} \biggr) \biggr),
   \quad k = 1\colon n.
\end{equation}
The \evals\ of $A_n$ are the reciprocals of the $\mu_k$. 

\subsection{The Kac--Murdock--Szeg\"o Matrix}\label{sec.KMS}

The \KMS\ matrix is the symmetric Toeplitz matrix,
depending on a single parameter $\rho\in\R$,
\begin{equation}\label{kms1}
A_n(\rho) = 
\begin{bmatrix}
   1          & \rho       & \rho^2 & \dots  & \rho^{n-1} \\[3pt]
   \rho       & 1          & \rho   & \dots  & \rho^{n-2} \\[-2pt]
   \rho^2     & \rho       & 1      & \ddots & \vdots     \\[-3pt]      
   \vdots     & \vdots     & \ddots & \ddots & \rho       \\[2pt]
   \rho^{n-1} & \rho^{n-2} & \dots  & \rho   & 1
\end{bmatrix} \in\mathbb{R}^{n\times n}.
\end{equation}
It was considered by Kac, Murdock, and Szeg\"o \cite[p.~784 ff.]{kms53},
who investigated its spectral properties.
It arises in the autoregressive AR(1) model in statistics and signal
processing.

It is straightforward to verify 
that
$A_n$ has a factorization $A_n = LDL^T$
with 
\begin{equation}\label{kms3}
 L = T_n(-\rho)^{-T}, \qquad
 D = \diag(1, 1-\rho^2, 1-\rho^2, \dots, 1-\rho^2).
\end{equation}
This \fact\ reveals several properties.
\begin{enumerate}[label=\upshape(\arabic*),wide,labelwidth=!]
\item 
$\det(A_n(\rho)) = (1-\rho^2)^{n-1}$.
\item
For $\rho \ne \pm 1$, $A_n$ is nonsingular and 
$A_n(\rho)^{-1} = T_n(-\rho) D^{-1} T_n(-\rho)^T$
is the tridiagonal (but not Toeplitz) matrix 
\begin{equation}\label{kms2}
  A_n(\rho)^{-1} = \frac{1}{1-\rho^2}
  \begin{bmatrix}
   1     & -\rho    &          &        &          &       \\[3pt]
   -\rho & 1+\rho^2 & -\rho    &        &          &       \\[-2pt]
         & -\rho    & 1+\rho^2 & \ddots &          &       \\[-2pt]      
         &          &  \ddots  & \ddots & \ddots   &       \\[2pt]   
         &          &          & -\rho  & 1+\rho^2 & -\rho \\[2pt]
         &          &          &        & -\rho    & 1
  \end{bmatrix}.
\end{equation}
\item
For $0 \le \rho \le 1$, 
$T_n(-\rho) = M(T_n(-\rho))$ and so by Theorem~\ref{thm.MTinv-TN}
$M(T_n(-\rho))^{-1} = T_n(-\rho)^{-1} = L^T$ is \TN,
so $A_n(\rho)$ is the product of three \TN\ matrices and is 
therefore \TN.
For $0< \rho < 1$,
$A_n(\rho)$ is also nonsingular and irreducible, 
so the \evals\ are distinct by Theorem~\ref{thm.TN-eig-irred-sing}.
Since $A_n(\rho) = \Sig A_n(-\rho) \Sig$ for $\Sig$ in \eqref{Sig-def}, 
$A_n(\rho)$ is similar to 
   $A_n(-\rho)$ and therefore
   $A_n(\rho)$ has distinct \evals\ for $0\ne \rho\in(-1,1)$.
\end{enumerate}

\subsection{The Pascal Matrix}\label{sec.pascal}

The Pascal matrix $P_n\in\R^{\nbyn}$, defined in \eqref{pascal-mat},
contains the rows of Pascal's triangle along the antidiagonals.
For example:
\begin{equation}\notag
 P_5 = 
\left[\begin{array}{ccccc} 
    1 & 1 & 1 & 1 & 1\\ 
    1 & 2 & 3 & 4 & 5\\ 
    1 & 3 & 6 & 10 & 15\\ 
    1 & 4 & 10 & 20 & 35\\ 
    1 & 5 & 15 & 35 & 70 \end{array}\right].
\end{equation}
This matrix is much-studied and most analyses involve the use of
combinatorial identities.
A number of key properties can be obtained from a \fact\ of $P_n$ into a
product of bidiagonal matrices.

The key observation is that $P_n$ can be reduced to upper \tri\ form by
repeatedly subtracting a  row from the row below.
For $n=5$, 
with $L_k(-1)$ denoting the unit lower bidiagonal matrix with $-1$s in
subdiagonal elements $(k+1,k)$, \dots, $(n-1,n)$,
\begin{align*}
   L_4(-1)L_3(-1)L_2(-1)L_1(-1) P_5 &= 
    \left[\begin{array}{ccccc} 1 & 0 & 0 & 0 & 0\\ 0 & 1 & 0 & 0 & 0\\ 0 & 0 & 1 & 0 & 0\\ 0 & 0 & 0 & 1 & 0\\ 0 & 0 & 0 & -1 & 1 \end{array}\right]
    \left[\begin{array}{ccccc} 1 & 0 & 0 & 0 & 0\\ 0 & 1 & 0 & 0 & 0\\ 0 & 0 & 1 & 0 & 0\\ 0 & 0 & -1 & 1 & 0\\ 0 & 0 & 0 & -1 & 1 \end{array}\right]\\
       & \hspace*{8pt} \times
     \left[\begin{array}{ccccc} 1 & 0 & 0 & 0 & 0\\ 0 & 1 & 0 & 0 & 0\\ 0 & -1 & 1 & 0 & 0\\ 0 & 0 & -1 & 1 & 0\\ 0 & 0 & 0 & -1 & 1 \end{array}\right]
    \left[\begin{array}{ccccc} 1 & 0 & 0 & 0 & 0\\ -1 & 1 & 0 & 0 & 0\\ 0 & -1 & 1 & 0 & 0\\ 0 & 0 & -1 & 1 & 0\\ 0 & 0 & 0 & -1 & 1 \end{array}\right]P_5\\
    &= 
    \left[\begin{array}{ccccc} 1 & 1 & 1 & 1 & 1\\
                               0 & 1 & 2 & 3 & 4\\ 
                               0 & 0 & 1 & 3 & 6\\ 
                               0 & 0 & 0 & 1 & 4\\ 
                               0 & 0 & 0 & 0 & 1 \end{array}\right]
       = R.
\end{align*}
In general, we have 
\begin{equation}\notag
     P_n = L_1(-1)^{-1} L_2(-1)^{-1} \dots L_{n-1}(-1)^{-1}R_n
         = L_n R_n,
\end{equation}
where $L_n$ is unit lower \tri\ and $R_n$ is unit upper \tri.
By the uniqueness of the LU and Cholesky \fact s of a positive definite matrix
we must have $L_n = R_n^T$,
so $P_n = R^T_nR_n$,
and it can be shown that
$R_n = L_{n-1}(1)^T L_{n-2}(1)^T\dots L_1(1)^T$,
which  contains the binomial coefficients downs its columns.

This is the \fact\ \eqref{TN-ALUfact} in Theorem~\ref{thm.TN-ALUfact}:
all the parameters are equal to $1$.

We can make several deductions.
\begin{enumerate}[label=\upshape(\arabic*),wide,labelwidth=!]

\item 
$P_n$ is \spd.

\item
$\det(P_n) = 1$.

\item
$P_n$ and $R_n$ are both \TN, since they are products of bidiagonal matrices
$L_i(1)$, each of which is \TN\ by Theorem~\ref{thm.tn-bidiag}.
Hence the \evals\ of $P_n$ are distinct by
Theorem~\ref{thm.TN-eig-irred-sing}.

\item
The matrix $S_n =\Sig R_n$ 
(where $\Sig$ is defined in \eqref{Sig-def})
is involutory, that is, $S_n^2 = I$.
This can be proved with the aid of the bidiagonal \fact\, but we omit the
rather tedious details.
Since $P_n  = S_n^TS_n$, we have 
$P_n^{-1} = S_n^{-1} S_n^{-T} = S_nS_n^T = S_n^{-T}P_n S_n^T$,
so $P_n^{-1}$ is similar to $P_n$, which means that the \evals\ of $P_n$
occur in reciprocal pairs.
It follows, in particular, that $\normt{P_n} = \normt{P_n^{-1}}$
and so $\ktwo(P_n) = \normt{P_n}^2$.

\end{enumerate}

It is also interesting to note that,
as an instance of Theorem~\ref{thm.exp-TNB},
the Cholesky factor $R_n$ is the
exponential of a bidiagonal matrix:
$R_n = \eu^{C_n}$, where \cite{actr01}, \cite{edst04}
\begin{equation}\notag
    C_n = 
     \begin{bmatrix} 
             0 & 1      &        &        &    \\
               & 0      & 2      &        &    \\
               &        & \ddots & \ddots &    \\
               &        &        & 0      & n-1 \\
               &        &        &        & 0
     \end{bmatrix} \in\R^{\nbyn}.
\end{equation}
The matrix $C_n$ is called the creation matrix
in \cite{acca17}, \cite{actr01} because of its role in generating 
matrix representations of \py s and providing simple proofs of identities.

\subsection{Tridiagonal Matrices from Partial Differential Equations}

Consider a linear system $Ax = b$, where $A = D+L+U$ with $D  = \diag(A)$
and $L$ and $U$ the
strictly lower \tri\ and 
strictly upper \tri\ parts of $A$, \resp.
The powers of the 
matrix $B = -(D+L)^{-1}U$
govern the convergence of the Gauss--Seidel iteration.
Note that $B$ is nonsymmetric and so in general can have complex \evals.

Suppose $A$ is tridiagonal with negative diagonal elements and nonnegative 
elements on the superdiagonal and subdiagonal, as is frequently the case 
in discretizations of partial differential equations, 
in which $A$ is typically a Toeplitz matrix.
For example, 
\begin{equation}\notag
   A = 
       \Bmatrix{-2 & 1  &    & \cr
                 1 & -2 &  1 & \cr
                   &  1 & -2 &1 \cr
                   &    &  1 & -2 }
\quad\Rightarrow\quad
   B = 
       \Bmatrix{
          0&   1/2&    0&    0\cr
          0&   1/4&  1/2&    0\cr
          0&   1/8&  1/4&  1/2\cr
          0&  1/16&  1/8&  1/4\cr
          }.
\end{equation}
The matrix $(-D -L)^{-1}$
is \TN\ by Theorem~\ref{thm.MTinv-TN}, because 
$-D-L = M(-D-L)$, 
and $U$ is \TN\ by Theorem~\ref{thm.tn-bidiag}.
Hence $B = (-D -L)^{-1} U$ is lower Hessenberg and \TN.
Furthermore, $B$ is irreducible if the subdiagonal of $L$ and the
superdiagonal of $U$ are nonzero.
Then Theorem~\ref{thm.TN-eig-irred-sing}
shows that the \evals\ of $B$ are real and nonnegative and the positive
\evals\ are distinct.
The \evals\ of $B$ can be deduced from the analysis 
of Young \cite{youn54}, \cite[Chap.~5]{youn71}.

\section*{Acknowledgement}

I thank Massimiliano Fasi and Xiaobo Liu
for their helpful comments on a draft manuscript.

\vspace*{-0.65cm}

\bibliographystyle{myplain2-doi}
\bibliography{strings,njhigham,la}

\end{document}